\RequirePackage{fix-cm}
\documentclass[smallextended]{svjour3}       
\smartqed  
\usepackage{amsmath, amscd, amsfonts, amssymb, graphicx, color}
\begin{document}

\title{Conditional $h$-convexity with applications  
}


\author{Ismail Nikoufar and Davuod Saeedi}


\institute{I. Nikoufar (Corresponding Author) \at
              Department of Mathematics, Payame Noor University, Tehran, Iran \\
              \email{nikoufar@pnu.ac.ir}           
\\
D. Saeedi \at
              Department of Mathematics, Payame Noor University, Tehran, Iran \\
              \email{dsaeedi3961@gmail.com}           
}

\date{Received: date / Accepted: date}

\maketitle

\begin{abstract}
In this paper, we introduce the notion of conditional $h$-convex functions and
we prove an operator version of the Jensen inequality for conditional $h$-convex functions.
Using this type of functions, we give some refinements for Ky-Fan's inequality, arithmetic-geometric mean inequality,
Chrystal inequality, and H$\ddot{o}$lder-McCarthy inequality.
Many of the other inequalities can be refined by applying this new notion.

\keywords{$h$-convex function\and Jensen's inequality\and Ky-Fan's inequality\and arithmetic-geometric mean inequality\and H\"{o}lder-McCarthy inequality.}

\subclass{26D15\and 15A39\and 47A63\and 46N10\and 47A60.}
\end{abstract}

\section{Introduction}

Throughout this paper, we denote by
$\mathcal{H}$ a Hilbert space and by $B(\mathcal{H})$ the algebra of all bounded linear operators on
$\mathcal{H}$.
The subalgebra of all self-adjoint operators in $B(\mathcal{H})$ is denoted by $B_{sa}(\mathcal{H})$.
An operator $A$ in $B_{sa}(\mathcal{H})$ is positive whenever $\langle Ax, x\rangle \geq 0$
for all $x\in\mathcal{H}$ and we write $A \geq 0$.
W denote by $Sp(A)$ the spectrum of an operator $A\in B(\mathcal{H})$.

One of the most important and most frequently applied in a
diversity of mathematical fields, especially in mathematical analysis and statistics
is Jensen's inequality.  In this process, it has been improved, generalized and adjusted
to various environments.
It has been linked to the other important inequalities and often
referred to as ``king of inequalities''.

The convexity of functions is an important issue in many fields of science, for instance in economy and optimization.
A function $f:\Bbb{I} \to \Bbb{R}$, $\Bbb{I}\subseteq \Bbb{R}$ is convex whenever the following inequality
$$
f(\lambda u + (1 - \lambda)v)\leq \lambda f(u) + (1- \lambda)f(v)
$$
holds for all $u, v \in \Bbb{I}$ and for all $\lambda\in [0,1]$ and the function $f : \Bbb{I} \to\Bbb{R}$ is concave whenever $-f$ is convex.

In 2007, Varo$\check{s}$anec \cite{S.Var-2007} introduced a wide class of functions
the so called $h$-convex functions which generalizes
convex, $s$-convex, Godunova-Levin, and $P$-class functions
\cite{Hudzik-AEqMath1994},
\cite{Dragomir-book2002,Nikoufar-CAOT2021,Pinheiro-AEqMath2007,Nikoufar-arXiv:2111.12837},
\cite{Godunova1985}, \cite{Dragomir-SJM1995}, \cite{Pearce1999}, \cite{Dragomir-BAMS1998,Nikoufar-Filomat2020}.

A non-negative function $f:\Bbb{I}\to\Bbb{R}$ is an $h$-convex function on $\Bbb{I}$ if
\begin{equation}\label{pi-eq}
f(\lambda u+(1-\lambda)v)\leq h(\lambda)f(u)+h(1-\lambda)f(v),
\end{equation}
where $h$ is a non-negative function defined on the real interval $\Bbb{J}$,
$u,v\in\Bbb{I}$ and $\lambda\in [0,1]\subseteq \Bbb{J}$.
For more results and generalizations regarding $h$-convexity, we refer the readers to see
\cite{Bombardelli2009,Dragomir-MM2015,Hazy2011,Olbrys2015,Nikoufar-arXiv:2201.05862}.


Jensen's inequality for convex  functions  is  one  of  the  most  important
result in the theory of inequalities and for appropriate choices of the function
many other famous inequalities are particular cases of this inequality.
An operator version of
the Jensen inequality for a convex function has been proved
by Mond and Pe$\check{c}$ari$\acute{c}$ as follows (\cite{Mond-HJM1993}, \cite{Furuta-Zagreb2005}):

\begin{theorem}\label{thm0}
Let $f:[m, M] \to \Bbb{R}$ be a continuous convex function.
Then,
\begin{equation}\label{eq0-thm1}
f(\langle Ax, x\rangle) \leq \langle f(A)x,x\rangle
\end{equation}
for every self-adjoint operator $A$ with $Sp(A)\subseteq [m,M]$ and every unit vector $x\in \mathcal{H}$.
\end{theorem}

In Section 2, we introduce the notion of conditional $h$-convex functions and
we prove an operator version of the Jensen inequality for conditional $h$-convex functions.
In Section 3,
Using this type of functions, we give some refinements for Ky-Fan's inequality, arithmetic-geometric mean inequality,
Chrystal inequality, and H$\ddot{o}$lder-McCarthy inequality.
We claim that many of the other inequalities can be refined by applying the notion of conditional $h$-convex functions.

%

\section{Mond and Pe$\check{c}$ari$\acute{c}$ inequality for conditional $h$-convex functions}



In this section, we provide our main results and
we indicate that an operator version of the Jensen inequality for conditional $h$-convex functions still holds
similar to that Mond and Pe$\check{c}$ari$\acute{c}$ considered for convex functions.

\begin{definition}
A function $f:[m,M]\to\Bbb{R}$ is conditional convex at $v\in [m,M]$ if
there exists a function $g:[m,M]\to\Bbb{R}$ such that $g(v)\leq v$ and
\begin{equation}\label{pi-eq}
f(\lambda u+(1-\lambda)v)\leq \lambda f(u)+(1-\lambda)f(v)
\end{equation}
for every $u\in [m,M]$ with $g(v)\leq u \leq v$ and $\lambda\in [0,1]$.
\end{definition}
It is clear that every convex function is conditional convex
at every $v\in [m,M]$ by considering $g(t)=m$ or $g(t)=\frac{m}{2}$, but not vice versa.

\begin{example}
Define $f:[0,2]\to\Bbb{R}$ by $f(t)=(t-1)^3$ and
$g(t)=\Big\{
\begin{array}{ll}
1, & t=2\\
2, & t\neq 2
\end{array}
$.
The function $f$ is not convex on $[0,2]$, where
$f$ is conditional convex at $2\in [0,2]$.
Consider $g(t)=\cos(\frac{4\pi}{3}t)$. Then $f$ is conditional convex at $\frac{3}{2}\in [0,2]$.
\end{example}

\begin{definition}
A non-negative function $f:[m,M]\to\Bbb{R}$ is conditional $h$-convex at $v\in [m,M]$ if
there exists a function $g:[m,M]\to\Bbb{R}$ such that $g(v)\leq v$ and
\begin{equation}\label{pi-eq}
f(\lambda u+(1-\lambda)v)\leq h(\lambda)f(u)+h(1-\lambda)f(v),
\end{equation}
for every $u\in [m,M]$ with $g(v)\leq u \leq v$ and $\lambda\in [0,1]$,
where $h$ is a non-negative function defined on the real interval $\Bbb{J}$, $[0,1]\subseteq \Bbb{J}$.
In this situation, we say the function $f$ is $g$-conditional $h$-convex at $v\in [m,M]$.
\end{definition}

It is clear that every $h$-convex function is $g$-conditional $h$-convex
at every $v\in [m,M]$ by considering $g(t)=m$, but not vice versa.

For a non-negative and non-zero function $h:K\subset \Bbb{R}\to\Bbb{R}$, we define
$$
M_K(h)=\inf_{t\in K}\frac{h(t)}{t}
$$
and call it the Jensen's coefficient for $h$-convex functions on $K$.
We consider all functions $h$ such that the Jensen's coefficient $M_K(h)$ exists.

\begin{theorem} \label{cnd-thm1}
Let $h:[0,1]\to\Bbb{R}$ be a non-negative and non-zero function.
If $f$ is a continuous $g$-conditional $h$-convex function at $v\in[m, M]$, then
\begin{equation}\label{eq0-thm1}
f(\langle Ax, x\rangle) \leq M_{(0,1)}(h)\langle f(A)x,x\rangle
\end{equation}
for each $x\in \mathcal{H}$ with $||x||=1$ and every self-adjoint operator $A$ such that $Sp(A)\subseteq [g(v),v]$.
\end{theorem}
\begin{proof}
It is clear that if $M_{(0,1)}(h)=+\infty$, then the inequality \eqref{eq0-thm1} holds.
Assume that $M_{(0,1)}(h)<+\infty$.
Since $f$ is $g$-conditional $h$-convex at $v\in [m, M]$,
\begin{equation}\label{eq1-thm1}
f(\lambda u+(1-\lambda)v)-h(1-\lambda)f(v)\leq h(\lambda) f(u)
\end{equation}
for all $u\in[g(v),v]$ and for all $\lambda\in (0,1)$.
By dividing both sides of \eqref{eq1-thm1} with $\lambda\in (0,1)$, one can reach
\begin{equation}\label{eq1-eq1-thm1}
\frac{f(\lambda u+(1-\lambda)v)-h(1-\lambda)f(v)}{\lambda}\leq \frac{h(\lambda)}{\lambda} f(u)
\end{equation}
for all $u\in[g(v),v]$.
Define
\begin{equation}\label{eq2-thm1}
\alpha:=\min_{u\in[g(v),v]}\frac{f(\lambda u+(1-\lambda)v)-h(1-\lambda)f(v)}{\lambda(u-v)}.
\end{equation}
The inequalities \eqref{eq1-eq1-thm1} and \eqref{eq2-thm1} entail that
$$
\alpha(u-v)\leq \frac{h(\lambda)}{\lambda}f(u)
$$
for all $u\in[g(v),v]$ and $\lambda\in (0,1)$.
Consider the linear function $l(t):=\alpha(t-v)$ and $\bar{g}=\langle Ax,x\rangle$
for every self-adjoint operator $A$ with $Sp(A)\subseteq [g(v),v]$.
This implies that $l(u)\leq \frac{h(\lambda)}{\lambda}f(u)$ for all $u\in[g(v),v]$ and
$g(v)\leq \bar{g}\leq v$, respectively.
Let $l'(t):=\alpha(t-\bar{g})+f(\bar{g})$ be the straight line passing through the point
$(\bar{g},f(\bar{g}))$ and parallel to the line $l$.
It follows from the continuity of the function $f$ that
\begin{equation}\label{eq3-thm1}
l'(\bar{g})\geq f(\bar{g})-\epsilon
\end{equation}
for all $\epsilon > 0$.
We now consider two cases:

(i) Assume that $l'(t)\leq \frac{h(\lambda)}{\lambda}f(t)$ for every $t\in [g(v),v]$.
By using the functional calculus, one has
$l'(A)\leq \frac{h(\lambda)}{\lambda}f(A)$
for every self-adjoint operator $A$ with $Sp(A)\subseteq [g(v),v]$ and so
\begin{equation}\label{eq4-thm1}
\langle l'(A)x,x \rangle\leq \frac{h(\lambda)}{\lambda}\langle f(A)x,x \rangle
\end{equation}
for each $x\in \mathcal{H}$ with $||x||=1$.
Since $l'$ is linear and $\epsilon$ is arbitrary the inequalities \eqref{eq3-thm1} and \eqref{eq4-thm1} imply
\begin{equation}\label{eq5-thm1}
f(\langle Ax,x\rangle)\leq \frac{h(\lambda)}{\lambda}\langle f(A)x,x \rangle.
\end{equation}

(ii) Assume that there exist some points $t\in[g(v),v]$ such that
$$
l'(t)>\frac{h(\lambda)}{\lambda}f(t).
$$
Define the sets $T$ and $S$ as follows:
\begin{align*}
T&:=\bigg\{t\in [g(v),\bar{g}]: l'(t)>\frac{h(\lambda)}{\lambda}f(t)\bigg\},\\
S&:=\bigg\{t\in [\bar{g},v]: l'(t)>\frac{h(\lambda)}{\lambda}f(t)\bigg\}.
\end{align*}
Consider $t_T:=\max\{t: t\in T\}$ and $t_S:=\min\{t: t\in S\}$.
Let $l_T$ be the line passing through the points $(t_T,0)$ and $(\bar{g},f(\bar{g}))$
and $l_S$ the line passing through the points $(t_S,0)$ and $(\bar{g},f(\bar{g}))$.
Define the function $L$ on $[g(v),v]$ as follows: 
$$
L(t):=\Big\{
\begin{array}{ll}
l_T(t), t\in [g(v),\bar{g}],\\
l_S(t), t\in [\bar{g},v].
\end{array}
$$
To show that $L(t)\leq \frac{h(\lambda)}{\lambda}f(t)$ holds for all $t\in [g(v), v]$,
consider the partition $\{g(v), t_T, \bar{g}, t_S, v\}$ for the closed interval $[g(v), v]$.
Note that $l_T(t)\leq 0$ for every $t\in[g(v), t_T]$ and
since $f(t)\geq 0$, we clearly observe that
$l_T(t)\leq \frac{h(\lambda)}{\lambda}f(t)$ for every $t\in[g(v), t_T]$.
On the other hand, we see that
\begin{equation}\label{main1}
l'(t)\leq \frac{h(\lambda)}{\lambda}f(t)
\end{equation}
for every $t\in (t_T, \bar{g}]$,
since if there exists $t_0\in (t_T, \bar{g}]$ such that
$l'(t_0)> \frac{h(\lambda)}{\lambda}f(t_0)$, then $t_0\in T$ and $t_0< t_T$, which is a contradiction.
So, by letting $t$ tends to $t_T$ from right in \eqref{main1}, the inequality \eqref{main1} indicates that
\begin{equation}\label{main110}
l'(t_T)\leq \frac{h(\lambda)}{\lambda}f(t_T).
\end{equation}
Moreover, since $t_T\in \bar{T}$, the reverse inequality holds in \eqref{main110} and hence
$l'$ is the line passing through the points $(t_T,\frac{h(\lambda)}{\lambda}f(t_T))$ and $(\bar{g},f(\bar{g}))$
and its slope is $\alpha'=\frac{f(\bar{g})-\frac{h(\lambda)}{\lambda}f(t_T)}{\bar{g}-t_T}$, where
the slope of $l_T$ is $\alpha_T=\frac{f(\bar{g})}{\bar{g}-t_T}$ and so $l_T(t)\leq l'(t)$ for every $t\in (t_T, \bar{g}]$.
By the inequality \eqref{main1} we observe that
$
l_T(t)
\leq \frac{h(\lambda)}{\lambda}f(t)
$
for every $t\in (t_T, \bar{g}]$.
So, $L(t)=l_T(t)\leq\frac{h(\lambda)}{\lambda}f(t)$ for every $t\in [g(v), \bar{g}]$.

By the similar methods one can show that $L(t)=l_S(t)\leq\frac{h(\lambda)}{\lambda}f(t)$ for every $t\in [\bar{g}, v]$.
Note that
the lines $l_T$ and $l_S$ are joining at the point along the length of $\bar{g}$ and so
$l_T(\bar{g})=l_S(\bar{g})$ and since $f$ is continuous,
\begin{equation}
l_T(\bar{g})=f(\bar{g})\geq f(\bar{g})-\epsilon
\end{equation}
for arbitrary $\epsilon > 0$. For the case  $Sp(A)\subseteq [g(v),\bar{g}]$, we have
\begin{align*}
f(\langle Ax,x\rangle)-\epsilon\leq l_T(\langle Ax,x\rangle)=\langle l_T(A)x,x\rangle\leq \frac{h(\lambda)}{\lambda}\langle f(A)x,x \rangle.
\end{align*}
Moreover, for the case $Sp(A)\subseteq [\bar{g},v]$, we have
\begin{align*}
f(\langle Ax,x\rangle)-\epsilon
&\leq l_T(\langle Ax,x\rangle)=l_S(\langle Ax,x\rangle)=\langle l_S(A)x,x\rangle\\
&\leq \frac{h(\lambda)}{\lambda}\langle f(A)x,x \rangle
\end{align*}
and consequently we can deduce \eqref{eq5-thm1}.
Now, by taking the infimum over all $\lambda\in (0,1)$ on both sides of the inequality \eqref{eq5-thm1},
one can deduce \eqref{eq0-thm1}.
\end{proof}

\begin{corollary}\label{thebest}
Under the hypotheses of Theorem \ref{cnd-thm1},
if the function $\frac{h(\lambda)}{\lambda}$ is decreasing on $(0,1)$, then $\lambda=\frac{1}{2}$
is the best possible in the inequality \eqref{eq5-thm1}. Indeed, one has
\begin{equation}\label{rem-eq}
f(\langle Ax, x\rangle) \leq 2h\bigg(\frac{1}{2}\bigg)\langle f(A)x,x\rangle
\end{equation}
for each $x\in \mathcal{H}$ with $||x||=1$ and every self-adjoint operator $A$ such that $Sp(A)\subseteq [g(v),v]$.
\end{corollary}
\begin{proof}
We divide the interval $(0,1)$ to two parts $(0,\frac{1}{2}]$ and $(\frac{1}{2},1)$, respectively.
If $\lambda\in (0,\frac{1}{2}]$, then
the infimum value of the function $\frac{h(\lambda)}{\lambda}$ occurs at the endpoint $\lambda=\frac{1}{2}$,
since the function $\frac{h(t)}{t}$ is decreasing and so
$$
M_{(0,\frac{1}{2}]}(h)=\frac{h\bigg(\frac{1}{2}\bigg)}{\frac{1}{2}}=2h\bigg(\frac{1}{2}\bigg).
$$
When $\lambda\in(\frac{1}{2},1)$,
we show that there is an $h$-convex function such that does not satisfy the inequality \eqref{eq5-thm1}.
Let $h(\lambda)=\lambda^{\beta}$, $0<\beta<1$, and $0<\lambda<1$.
Then, $\frac{h(\lambda)}{\lambda}$ is decreasing on $(0,1)$.
Define $f: [0, \infty)\to \Bbb{R}$ by $f(t)=e^{-t}$.
The function $f$ is $h$-convex, since by using the convexity of $f$, we have
\begin{align*}
f(\lambda x+(1-\lambda)y)&=e^{-(\lambda x+(1-\lambda)y)}\\
&\leq \lambda e^{-x}+(1-\lambda)e^{-y}\\
&\leq\lambda^{\beta} e^{-x}+(1-\lambda)^{\beta}e^{-y}\\
&=h(\lambda)f(x)+h(1-\lambda)f(y)
\end{align*}
for every $x,y\geq 0$ and $\lambda\in [0,1]$.
Consider
$
A=\Big(
\begin{array}{ll}
0 & 0\\
0 & a
\end{array}
 \Big)
$, $a>0$,
and $x=(\frac{1}{\sqrt{2}},\frac{1}{\sqrt{2}})$.
A simple calculation shows that $f(\langle Ax, x\rangle)=f(\frac{a}{2})=e^{-\frac{a}{2}}$ and $\langle f(A)x,x\rangle=\frac{e^{-a}}{2}$.
By applying the inequality \eqref{eq5-thm1}, one can see that
$e^{-\frac{a}{2}}\leq\lambda^{\beta-1}\frac{e^{-a}}{2}$
which implies $2<2e^{\frac{a}{2}}\leq\lambda^{\beta-1}$.
On the other hand, for $\lambda\in(\frac{1}{2},1)$, we have $\lambda^{\beta-1}<\lambda^{-1}<2$ which is a contradiction.
\end{proof}

\section{Some Applications}
In this section, we provide some examples of conditional $h$-convex functions.
Using this type of functions, we give some refinements for Ky-Fan's inequality, arithmetic-geometric mean inequality,
Chrystal inequality, and H$\ddot{o}$lder-McCarthy inequality.

%
%
%
\begin{lemma}\label{Ky-Fan}
Let $h:[0,1]\to (0,\infty)$ defined by $h(t)=\frac{\alpha}{\beta}e^{t(1-t)}$
and
let $g:(0,\frac{1}{2}]\to[0,\infty)$ defined by
$g(t)=\frac{t^{\alpha}}{t^{\alpha}+(1-t)^{\alpha}}$
with $1<\alpha\leq\beta\leq\alpha+1$.
Then, the function $f:(0,\frac{1}{2}]\to[0,\infty)$ defined by $f(t)=\ln(\frac{1-t}{t})$
is $g$-conditional $h$-convex at every $v\in(0,\frac{1}{2}]$.
\end{lemma}
\begin{proof}
Let $v\in(0,\frac{1}{2}]$ and define the function $F$ by
$$
F(u,\lambda)
=\frac{\alpha}{\beta}e^{t(1-t)}\ln\bigg(\frac{1-u}{u}\bigg)
+\frac{\alpha}{\beta}e^{t(1-t)}\ln\bigg(\frac{1-v}{v}\bigg)
-\ln\bigg(\frac{1-(\lambda u+(1-\lambda)v)}{\lambda u+(1-\lambda)v}\bigg)
$$
for every $u\in[g(v),v]$ and $\lambda\in[0,1]$.
It is enough to show that $F(u,\lambda)\geq 0$.
We have
\begin{align}\label{setare}
F(u,\lambda)
&\geq\frac{\alpha}{\beta}\ln\bigg(\frac{1-u}{u}\bigg)
    +\frac{\alpha}{\beta}\ln\bigg(\frac{1-v}{v}\bigg)
    -\ln\bigg(\frac{1-(\lambda u+(1-\lambda)v)}{\lambda u+(1-\lambda)v}\bigg)\nonumber\\
&=\ln\bigg(\frac{1-u}{u}\bigg)
 -\ln\bigg(\frac{1-(\lambda u+(1-\lambda)v)}{\lambda u+(1-\lambda)v}\bigg)\nonumber\\
&\ \ \ +\frac{\alpha}{\beta}\ln\bigg(\frac{1-v}{v}\bigg)
 -(1-\frac{\alpha}{\beta})\ln\bigg(\frac{1-u}{u}\bigg).
\end{align}
It is obvious that $u<\lambda u+(1-\lambda)v$ for every $\lambda\in[0,1]$.
So, $\frac{1-(\lambda u+(1-\lambda)v}{\lambda u+(1-\lambda)v}<\frac{1-u}{u}$.
It follows that
\begin{equation}\label{yek}
\ln\bigg(\frac{1-u}{u}\bigg)-\ln\bigg(\frac{1-(\lambda u+(1-\lambda)v)}{\lambda u+(1-\lambda)v}\bigg)\geq 0.
\end{equation}
On the other hand, we have $\frac{v^{\alpha}}{v^{\alpha}+(1-v)^{\alpha}}\leq u$ so that
$\frac{1-u}{u}\leq (\frac{1-v}{v})^{\alpha}$.
This ensures
$\ln(\frac{1-u}{u})\leq {\alpha}\ln(\frac{1-v}{v})$.
Since $\beta\leq \alpha+1$,
this indicates that
\begin{equation}\label{doo}
\bigg(1-\frac{\alpha}{\beta}\bigg)\ln\bigg(\frac{1-u}{u}\bigg)\leq \frac{1}{\beta}\ln\bigg(\frac{1-u}{u}\bigg)\leq \frac{\alpha}{\beta}\ln\bigg(\frac{1-v}{v}\bigg).
\end{equation}
The inequalities \eqref{yek} and \eqref{doo} entail that $F(u,\lambda)\geq 0$.
\end{proof}

\begin{corollary}\label{Ky-Fan-cor}
(Refined Ky-Fan's inequality)
If $\alpha>1$, $v\in (0,\frac{1}{2}]$, and
$a_i\in\bigg[\frac{v^{\alpha}}{v^{\alpha}+(1-v)^{\alpha}},v\bigg]$,
then
\begin{equation}\label{eq0-Ky-Fan}
\frac{\sum_{i=1}^{n}q_i(1-a_i)}{\sum_{i=1}^{n}q_ia_i}
\leq\prod_{i=1}^{n}\bigg(\frac{1-a_i}{a_i}\bigg)^{\frac{\alpha}{\alpha+\gamma}q_i}
\leq\prod_{i=1}^{n}\bigg(\frac{1-a_i}{a_i}\bigg)^{q_i},
\end{equation}
where $q_i\in (0,1)$, $\sum_{i=1}^{n}q_i=1$, and $\gamma=\max|a_i-a_j|$.
The equality occurs when $\gamma=0$.
\end{corollary}
\begin{proof}
Consider $f(t)=\ln(\frac{1-t}{t})$ on $(0,\frac{1}{2}]$.
Since $0<a_i\leq v\leq \frac{1}{2}$, we have $1<\alpha\leq\alpha+\gamma\leq\alpha+1$.
By Lemma \ref{Ky-Fan} and setting $\beta=\alpha+\gamma$, the function $f$ is $g$-conditional $h$-convex
at every $v\in(0,\frac{1}{2}]$.
Define
$$
A=\left(
\begin{array}{lll}
a_1 & \cdots & 0\\
\vdots & \ddots & \vdots\\
0 & \cdots & a_n
\end{array}
\right),\ \ \
x=\left(
\begin{array}{l}
\sqrt{q_1}\\
\vdots\\
\sqrt{q_n}
\end{array}
\right).
$$
Then, $Sp(A)\subseteq \bigg[\frac{v^{\alpha}}{v^{\alpha}+(1-v)^{\alpha}},v\bigg]$ and $||x||=1$.
According to Theorem \ref{cnd-thm1} one knows
\begin{equation}\label{eq1-Ky-Fan}
f(\langle Ax, x\rangle) \leq M_{(0,1)}(h)\langle f(A)x,x\rangle.
\end{equation}
On the other hand, we have
\begin{align*}
\langle Ax, x\rangle&=\sum_{i=1}^{n}q_ia_i,\\
M_{(0,1)}(h)&=\inf_{\lambda\in [0,1]}\frac{h(\lambda)}{\lambda}=\frac{\alpha}{\beta},\\
\langle f(A)x,x\rangle&=\sum_{i=1}^{n}q_i\ln\bigg(\frac{1-a_i}{a_i}\bigg).
\end{align*}
By replacing these quantities in \eqref{eq1-Ky-Fan}, we get
the first inequality in \eqref{eq0-Ky-Fan}.
The second inequality follows from the fact that
$\frac{1-a_i}{a_i}\geq 1$ and $\frac{\alpha}{\beta}\leq 1$.
Note that when $\gamma=0$, we have $a_i=a_j=a$ for every $i\neq j$ and the equality occurs.
\end{proof}

%
%
%

\begin{lemma}\label{Arth-Geo}
Let $h:[0,1]\to (0,\infty)$ defined by $h(t)=\frac{\alpha}{\beta}e^{t(1-t)}$
and
let $g:(0,1]\to[0,\infty)$ defined by $g(t)=t^{\alpha}$ with $1<\alpha\leq\beta\leq\alpha+1$.
Then, the function $f:(0,1]\to[0,\infty)$ defined by $f(t)=-\ln(t)$
is $g$-conditional $h$-convex at every $v\in(0,1]$.
\end{lemma}
\begin{proof}
Let $v\in(0,1]$ and define the function $F$ by
$$
F(u,\lambda)
=-\frac{\alpha}{\beta}e^{t(1-t)}\ln(u)
-\frac{\alpha}{\beta}e^{t(1-t)}\ln(v)
+\ln(\lambda u+(1-\lambda)v)
$$
for every $u\in[v^{\alpha},v]$ and $\lambda\in[0,1]$.
We demonstrate that $F(u,\lambda)\geq 0$.
We have
\begin{align*}
F(u,\lambda)
&\geq -\frac{\alpha}{\beta}\ln(u)
    -\frac{\alpha}{\beta}\ln(v)
    +\ln(\lambda u+(1-\lambda)v)\nonumber\\
&=\bigg(1-\frac{\alpha}{\beta}\bigg)\ln(u)-\frac{\alpha}{\beta}\ln(v)-\ln(u)+\ln(\lambda u+(1-\lambda)v).
\end{align*}
Since $u<\lambda u+(1-\lambda)v$ for every $\lambda\in[0,1]$, $\ln(u)<\ln(\lambda u+(1-\lambda)v)$
and hence
\begin{equation}\label{App2-1}
-\ln(u)+\ln(\lambda u+(1-\lambda)v)>0.
\end{equation}
On the other hand, we have $v^{\alpha}\leq u$.
Since $\beta-\alpha\leq 1$ and $0<u<1$, $u\leq u^{\beta-\alpha}$.
These two inequalities ensure that $v^{\alpha}\leq u^{\beta-\alpha}$.
So,
\begin{equation}\label{App2-2}
\frac{\alpha}{\beta}\ln(v)\leq \frac{\beta-\alpha}{\beta}\ln(u).
\end{equation}
We realize from the inequalities \eqref{App2-1} and \eqref{App2-2} that $F(u,\lambda)\geq 0$.
\end{proof}

\begin{corollary}\label{Arth-Geo-cor}
(Refined Arithmetic-Geometric Mean inequality)
If $\alpha>1$, $v\in (0,1]$, and $a_i\in [v^{\alpha},v]$,
then
\begin{equation}\label{eq0-Arth-Goe}
\prod_{i=1}^{n}a_i^{q_i} \leq \prod_{i=1}^{n}a_i^{\frac{\alpha}{\alpha+\gamma}q_i}\leq\sum_{i=1}^{n}q_ia_i,
\end{equation}
where $q_i\in (0,1)$, $\sum_{i=1}^{n}q_i=1$, and $\gamma=\max|a_i-a_j|$.
The equality occurs when $\gamma=0$.
\end{corollary}
\begin{proof}
Consider $f(t)=-\ln(t)$ on $(0,1]$.
Since $0<a_i\leq v\leq 1$, we have $1<\alpha\leq\alpha+\gamma\leq\alpha+1$. 
Due to Lemma \ref{Ky-Fan} and setting $\beta=\alpha+\gamma$, the function $f$ is $g$-conditional $h$-convex
at every $v\in(0,1]$.
Define the matrix $A$ and the vector $x$ as in the proof of Corollary \ref{Ky-Fan-cor}.
Then, $Sp(A)\subseteq [v^{\alpha},v]$ and $||x||=1$.
Theorem \ref{cnd-thm1} signifies the inequality \eqref{eq1-Ky-Fan} holds and by the assumptions we have
\begin{align*}
\langle Ax, x\rangle&=\sum_{i=1}^{n}q_ia_i,\\
M_{(0,1)}(h)&=\inf_{\lambda\in [0,1]}\frac{h(\lambda)}{\lambda}=\frac{\alpha}{\beta},\\
\langle f(A)x,x\rangle&=-\sum_{i=1}^{n}q_i\ln(a_i).
\end{align*}
By replacing these quantities in \eqref{eq1-Ky-Fan}, one can deduce
\begin{equation*}
-\ln\bigg(\sum_{i=1}^{n}q_ia_i\bigg)\leq -\frac{\alpha}{\beta}\sum_{i=1}^{n}q_i\ln(a_i).
\end{equation*}
From this inequality, it follows that
\begin{equation*}
\ln\bigg(\sum_{i=1}^{n}q_ia_i\bigg)\geq \sum_{i=1}^{n}\ln\bigg(a_i^{\frac{\alpha}{\beta}q_i}\bigg)=\ln\bigg(\prod_{i=1}^{n}a_i^{\frac{\alpha}{\beta}q_i}\bigg),
\end{equation*}
which is equivalent to the second inequality in \eqref{eq0-Arth-Goe}.
The first inequality comes from the fact that
$0<a_i\leq 1$ and $\frac{\alpha}{\beta}\leq 1$.
Note that when $\gamma=0$, we have $a_i=a_j=a$ for every $i\neq j$ and the equality occurs.
\end{proof}


%
%
%

\begin{lemma}\label{Chrystal}
Let $h:[0,1]\to (0,\infty)$ defined by $h(t)=\frac{\alpha}{\beta}e^{t(1-t)}$
and
let $g:(0,\infty)\to[0,\infty)$ defined by $g(t)=\ln((1+e^t)^{\frac{\beta}{\alpha}-1}-1)$ with $0<\alpha\leq\beta\leq 2\alpha$.
Then, the function $f:(0,\infty)\to[0,\infty)$ defined by $f(t)=\ln(1+e^t)$
is $g$-conditional $h$-convex at every $v\in(0,\infty)$.
\end{lemma}
\begin{proof}
Let $v\in(0,\infty)$ and define the function $F$ by
$$
F(u,\lambda)
=\frac{\alpha}{\beta}e^{t(1-t)}\ln(1+e^u)
+\frac{\alpha}{\beta}e^{t(1-t)}\ln(1+e^v)
-\ln(1+e^{\lambda u+(1-\lambda)v})
$$
for every $u\in\bigg[\ln((1+e^v)^{\frac{\beta}{\alpha}-1}-1),v\bigg]$ and $\lambda\in[0,1]$.
We prove that $F(u,\lambda)\geq 0$.
Note that
\begin{align*}
F(u,\lambda)
&\geq \frac{\alpha}{\beta}\ln(1+e^u)
+\frac{\alpha}{\beta}\ln(1+e^v)
-\ln(1+e^{\lambda u+(1-\lambda)v}).
\end{align*}
Since $\lambda u+(1-\lambda)v<v$ for every $\lambda\in[0,1]$,
\begin{equation}\label{App3-1}
\ln(1+e^{\lambda u+(1-\lambda)v})<\ln(1+e^v).
\end{equation}
Since $0<\frac{\beta}{\alpha}-1\leq 1$, $\ln((1+e^v)^{\frac{\beta}{\alpha}-1}-1)>0$.
On the other hand, we have $\ln((1+e^v)^{\frac{\beta}{\alpha}-1}-1)\leq u$. So,
\begin{equation}\label{App3-2}
1+e^{v}\leq (1+e^u)^{\frac{\alpha}{\beta}}(1+e^v)^{\frac{\alpha}{\beta}}.
\end{equation}
In view of the inequalities \eqref{App3-1} and \eqref{App3-2}, we find that
\begin{equation*}
\ln(1+e^{\lambda u+(1-\lambda)v})<\ln((1+e^u)^{\frac{\alpha}{\beta}}(1+e^v)^{\frac{\alpha}{\beta}}),
\end{equation*}
which implies $F(u,\lambda)\geq 0$.
\end{proof}

\begin{corollary}\label{Chrystal-cor}
(Refined Chrystal's inequality)
If $\alpha>0$, $0<v\leq\alpha$, and $a_i, b_i\in [\ln((1+e^v)^{\frac{\alpha+\gamma}{\alpha}-1}-1),v]$,
then
\begin{equation}\label{eq0-Chrystal}
\prod_{i=1}^{n}a_i^{q_i} + \prod_{i=1}^{n}b_i^{q_i}
\leq \prod_{i=1}^{n}\bigg(\frac{(a_i+b_i)^{\frac{\alpha}{\alpha+\gamma}}}{b_i^{\frac{\alpha}{\alpha+\gamma}-1}}\bigg)^{q_i}
\leq \prod_{i=1}^{n}(a_i+b_i)^{q_i},
\end{equation}
where $q_i\in (0,1)$, $\sum_{i=1}^{n}q_i=1$, and
$$\gamma=\max\bigg\{\max|a_i-a_j|, \max|b_i-b_j|, \max|a_i-b_j|\bigg\}.$$
The equality occurs when $\gamma=0$.
\end{corollary}
\begin{proof}
Consider $f(t)=\ln(1+e^t)$ on $(0,\infty)$.
Since $0<a_i\leq v\leq\alpha$, we have $0<\alpha\leq\alpha+\gamma\leq 2\alpha$.
So, the function $f$ is $g$-conditional $h$-convex
at every $v\in (0,\alpha]$, by Lemma \ref{Chrystal} and putting $\beta=\alpha+\gamma$.
Defining the matrix $A$ and the vector $x$ as in the proof of Corollary \ref{Ky-Fan-cor}
ensures $Sp(A)\subseteq [\ln((1+e^v)^{\frac{\beta}{\alpha}-1}-1),v]$ and $||x||=1$.
In view of Theorem \ref{cnd-thm1} we know that the inequality \eqref{eq1-Ky-Fan} holds and by the assumptions we have
\begin{align*}
\langle Ax, x\rangle&=\sum_{i=1}^{n}q_ia_i,\\
M_{(0,1)}(h)&=\inf_{\lambda\in [0,1]}\frac{h(\lambda)}{\lambda}=\frac{\alpha}{\beta},\\
\langle f(A)x,x\rangle&=\sum_{i=1}^{n}q_i\ln(1+e^{a_i}).
\end{align*}
By replacing these quantities in \eqref{eq1-Ky-Fan}, we see that
\begin{equation}\label{eq3-Chrystal}
\ln\bigg(1+e^{\sum_{i=1}^{n}q_ia_i}\bigg)\leq \frac{\alpha}{\beta}\sum_{i=1}^{n}q_i\ln(1+e^{a_i}).
\end{equation}
On the other hand, since $\frac{\alpha}{\beta}\leq 1$,
\begin{equation}\label{eq4-Chrystal}
\frac{\alpha}{\beta}\sum_{i=1}^{n}q_i\ln(1+e^{a_i})\leq \sum_{i=1}^{n}q_i\ln(1+e^{a_i}).
\end{equation}
Combining the inequalities \eqref{eq3-Chrystal} and \eqref{eq4-Chrystal}, we get
\begin{equation}\label{eq5-Chrystal}
\ln\bigg(1+e^{\sum_{i=1}^{n}q_ia_i}\bigg)\leq \frac{\alpha}{\beta}\sum_{i=1}^{n}q_i\ln(1+e^{a_i})
\leq\sum_{i=1}^{n}q_i\ln(1+e^{a_i}).
\end{equation}
By replacing $\ln(\frac{a_i}{b_i})$ with $a_i$ in \eqref{eq5-Chrystal}, we observe that
\begin{equation*}
\ln\bigg(1+\prod_{i=1}^{n}(\frac{a_i}{b_i})^{q_i}\bigg)
\leq \ln\bigg(\prod_{i=1}^{n}(1+\frac{a_i}{b_i})^{\frac{\alpha}{\beta}q_i}\bigg)
\leq\ln\bigg(\prod_{i=1}^{n}(1+\frac{a_i}{b_i})^{q_i}\bigg),
\end{equation*}
which implies
\begin{equation*}
1+\prod_{i=1}^{n}\bigg(\frac{a_i}{b_i}\bigg)^{q_i}
\leq \prod_{i=1}^{n}\bigg(1+\frac{a_i}{b_i}\bigg)^{\frac{\alpha}{\beta}q_i}
\leq\prod_{i=1}^{n}\bigg(1+\frac{a_i}{b_i}\bigg)^{q_i}.
\end{equation*}
Multiplying both sides by $\prod_{i=1}^{n}{b_i}^{q_i}$ gives \eqref{eq0-Chrystal}.
Note that when $\gamma=0$, we have $a_i=b_j=a$ for every $i\neq j$ and the equality occurs.
\end{proof}


%
%
%

\begin{lemma}\label{HM}
Let $h:[0,1]\to (0,\infty)$ defined by $h(t)=\frac{\alpha}{\beta}e^{t(1-t)}$
and
let $g:[0,\infty)\to[0,\infty)$ defined by $g(t)=t\sqrt[p]{\frac{\beta}{\alpha}-1}$
with $p>1$ and $0<\alpha\leq\beta\leq 2\alpha$.
Then, the function $f:[0,\infty)\to[0,\infty)$ defined by $f(t)=t^p$
is $g$-conditional $h$-convex at every $v\in [0,\infty)$.
\end{lemma}
\begin{proof}
Suppose $v\in [0,\infty)$.
Since $\beta\leq 2\alpha$, $v\sqrt[p]{\frac{\beta}{\alpha}-1}\leq v$.
Consider the function $F$ by
$$
F(u,\lambda)
=\frac{\alpha}{\beta}e^{t(1-t)}u^p
+\frac{\alpha}{\beta}e^{t(1-t)}v^p
-({\lambda u+(1-\lambda)v})^p
$$
for every $u\in\bigg[v\sqrt[p]{\frac{\beta}{\alpha}-1},v\bigg]$ and $\lambda\in[0,1]$.
To show that $F(u,\lambda)\geq 0$, we observe that
\begin{align*}
F(u,\lambda)
&\geq \frac{\alpha}{\beta}u^p+\frac{\alpha}{\beta}v^p-({\lambda u+(1-\lambda)v})^p\\
&=\frac{\alpha}{\beta}u^p-\bigg(1-\frac{\alpha}{\beta}\bigg)v^p+v^p-({\lambda u+(1-\lambda)v})^p.
\end{align*}
Since $\lambda u+(1-\lambda)v<v$ for every $\lambda\in[0,1]$,
\begin{equation}\label{App4-1}
(\lambda u+(1-\lambda)v)^p<v^p.
\end{equation}
Since 
$v\sqrt[p]{\frac{\beta}{\alpha}-1}\leq u$, $v^p(\beta-\alpha)\leq\alpha u^p$
and so dividing both sides by $\beta$ one has
\begin{equation}\label{App4-2}
v^p(1-\frac{\alpha}{\beta})\leq\frac{\alpha}{\beta} u^p.
\end{equation}
The inequalities \eqref{App4-1} and \eqref{App4-2} signify
$F(u,\lambda)\geq 0$.
\end{proof}

\begin{corollary}\label{HM-cor}
(Refined H$\ddot{o}$lder-McCarthy inequality)
If $\alpha>0$, $0<v\leq\alpha$, and $p>1$, then
\begin{equation*}
\langle Ax,x\rangle^p\leq \frac{\alpha}{\alpha+\gamma}\langle A^px,x\rangle \leq \langle A^px,x\rangle
\end{equation*}
for every self-adjoint operator $A$ with $Sp(A)\subseteq \bigg[v\sqrt[p]{\frac{\alpha+\gamma}{\alpha}-1},v\bigg]$
and every unit vector $x$, where $\gamma=\max\{|\eta-\zeta|: \eta,\zeta\in Sp(A)\}$.
The equality occurs when $\gamma=0$.
\end{corollary}
\begin{proof}
Define $f(t)=t^p$ on $(0,\infty)$.
Since $0<\eta\leq v\leq\alpha$ for every $\eta\in Sp(A)$, we have $0<\alpha\leq\alpha+\gamma\leq 2\alpha$.
Due to Lemma \ref{HM} and putting $\beta=\alpha+\gamma$,
the function $f$ is $g$-conditional $h$-convex
at every $v\in (0,\alpha]$.
In account of Theorem \ref{cnd-thm1} the inequality \eqref{eq1-Ky-Fan} holds and we have
\begin{align*}
\langle Ax, x\rangle^p \leq M_{(0,1)}(h)\langle A^px,x\rangle,
\end{align*}
where $M_{(0,1)}(h)=\inf_{\lambda\in [0,1]}\frac{h(\lambda)}{\lambda}=\frac{\alpha}{\beta}$.
Since $\frac{\alpha}{\beta}<1$, we get the desired result.
If $\gamma=0$, then $A$ is a multiple of the identity and one has the equality.
\end{proof}

\section{Declarations}
We remark that the potential conflicts of interest and
data sharing not applicable to this article and no data sets were generated during the current study.




\end{document}